\documentclass[11pt]{article}
\usepackage{graphics,latexsym,amssymb,amsmath,amscd}
\RequirePackage{graphics}
\usepackage{graphicx}
\usepackage{psfrag}
\input epsf
\date{}

\def\Im{\hbox{\rm Im}\,}
\def\[#1\]{\begin{eqnarray*}#1\end{eqnarray*}}

\def\ch#1{{{\bf H}^{#1}_{\C}}}

\def\Heis{{\mathfrak N}}

\def\phi{\varphi}

\def\Im{\hbox{\rm Im}\,}
\def\[#1\]{\begin{eqnarray*}#1\end{eqnarray*}}

\def\ch#1{{{\bf H}^{#1}_{\C}}}

\def\Heis{{\mathfrak N}}

\def\phi{\varphi}

\newtheorem{thm}{Theorem}[section]
\newtheorem{dfn}[thm]{Definition}

\newtheorem{prop}[thm]{Proposition}

\newcommand{\Pf}{{\sc Proof}. }
\newcommand{\EPf}{\hbox{}\hfill$\Box$\vspace{.5cm}}

\def\A{{\mathbb A}}

\input xy
\xyoption{all}

\newcommand{\C}{{{\mathbb C}}}
\newcommand{\R}{{{\mathbb R}}}

\newcommand{\PP}{{{\mathbb P}}}
\begin{document}
\title{A Volume function for Spherical CR tetrahedra}
\author{ Elisha Falbel\\
Institut de Math{\'e}matiques\\
Universit{\'e} Pierre et Marie Curie\\
4, place Jussieu \\ F-75252 Paris\\ 
e-mail: {\tt falbel\@math.jussieu.fr}}
\maketitle
\begin{abstract}
We define a volume function on configurations of
four points in the sphere $S^3$ which is invariant under 
the action of $PU(2,1)$, the automorphism group of the
CR structure defined on $S^3$ by its embedding in $\C^2$.
We show that the volume function, constructed using appropriate
combinations of the dilogathm function of Bloch-Wigner,
 satisfies a five term relation
in a  more general context which includes at the same time
CR  and real hyperbolic geometry.
\end{abstract}

\section{Introduction}

The volume of ideal tetrahedra in real three dimensional 
hyperbolic space is a fundamental geometric invariant.  One can interpret it as a function
defined  on ordered quadruples of distinct points 
in the Riemann sphere $\C P^1$.  It is invariant by
the diagonal action of $PSL(2,\C)$
 and  it satisfies a 5 term relation.
It defines an element of $H_{cont}^3(PSL(2,\C),\R)$ (the continuous group cohomology of $PSL(2,\C)$).

The goal of this paper is to
 define a volume function defined on ordered quadruples of points in a
more general context which includes at the same time real hyperbolic
geometry and CR geometry.  In particular it is defined for ordered quadruples of points
in the sphere $S^3\subset \C^2$.  It is invariant under the diagonal action
of $PU(2,1)$ on the configuration of points and satisfies a 5 term relation.
Although $H_{cont}^3(PU(2,1),\R)=0$ and therefore, in that case, this volume is a coboundary
 it is sufficiently interesting to be singled out.

Recall that ordered triples of (pairwise distinct) points
in $\C$ are classified, up to similarity, by a complex parameter
$z\in \C \setminus \{ 0,1\}$. A natural compactification
of this space, being obtained collapsing pairs of points, is identified
to $\C P^1$.

 More precisely, to each vertex of an ordered triple $p_1,p_2,p_3\in \C$
 we associate a coordinate: 
$z_1=\frac{p_3-p_1}{p_2-p_1},
z_2=\frac{p_1-p_2}{p_3-p_2},
z_3=\frac{p_2-p_3}{p_1-p_3}$.  They satisfy
$$
z_2=\frac{1}{1-z_1}, \ \ \ z_3=1-\frac{1}{z_1}.
$$ 
One can chose one of the coordinates (say $z_1\in  \C \setminus \{ 0,1\}$)
  to parametrize the configuration 
up to similarity. 

The main idea of this paper is that each point in a 
configuration of four points might see the other three points as 
forming an Euclidean triangle.   The mean value of the 
Bloch-Wigner dilogarithm of the invariants of
those four triangles satisfies a five term relation if natural compatibilities
are imposed between the triangles (up to similarity) associated to
a configuration of five points. 
 
 For a compact manifold the total 
volume is a function on an algebraic variety of ``geometric structures'' given by invariants $z_{ij}$ (satisfying certain compatibility conditions) associated to simplices of a triangulation. The volume at a hyperbolic structure  coincides with
the hyperbolic volume and the volume at closed spherical CR structure is always null.
I thank Riccardo Benedetti, Herbert Gangl, Juliette Genzmer, Julien March\'e, John Parker, Luc Pirio,  Qingxue Wang and Pierre Will
 for the discussions leading to that 
paper.
\section{Tetrahedra with cross-ratio structures}

Consider a set of four elements $\Delta= \{ p_1,p_2,p_3,p_4 \}$.  We call  $p_i$, $1\leq i\leq 4$ 
the vertices of $\Delta$.  Let $O\Delta$ be the set of all orderings of $\Delta$.  We will denote an element of $O\Delta$ by $[p_i,p_j,p_k,p_l]$ (where $\{ i,j,k,l\}=\{ 1,2,3,4\}$) and call it a simplex although we only deal with configurations of four points. Given $\Delta$,
there are 24 simplices divided in two classes $O\Delta^+$ (containing $[p_1,p_2,p_3,p_4]$)
 and $\Delta^-$ (containing $[p_1,p_2,p_4,p_3]$)
  of 12 elements each.  Each class is an orbit of the even permutation group acting on $O\Delta$.  

The following definition assigns similarity invariants
to each vertex of a configuration of four points.  That is, to each vertex, it is assigned
a triangle in $\C$,  up to similarity, defined by coordinates as in the introduction.

\begin{dfn} A cross-ratio structure on a set of four points  $\Delta=\{ p_1,p_2,p_3,p_4 \}$ is a function defined on the ordered quadruples 
$$
{\bf X} : O\Delta\rightarrow \C\setminus \{0,1\}
$$
satisfying, if 
$(i,j,k,l)$ is any permutation of $(1,2,3,4)$, the relations 
\begin{enumerate}
\item 
$$
{\bf X}(p_i,p_j,p_k,p_l)=\frac{1}{{\bf X}(p_i,p_j,p_l,p_k)}.
$$
\item  (similarity relations)
$$
{\bf X}(p_i,p_j,p_k,p_l)=\frac{1}{1-{\bf X}(p_i,p_l,p_j,p_k)}.
$$
\end{enumerate}
\end{dfn}

\begin{figure}
\setlength{\unitlength}{1cm}
\begin{center}
\begin{picture}(7,10)
\psfrag{z1}{$z_{12}$}
\psfrag{z1'}{$z_{21}$}
\psfrag{z~1}{$z_{34}$}
\psfrag{z~1'}{$z_{43}$}
\psfrag{z2}{$z_{13}$}
\psfrag{z2'}{$z_{24}$}
\psfrag{z~2}{$ z_{31}$}
\psfrag{z~2'}{$ z_{42}$}
\psfrag{z3}{$z_{14}$}
\psfrag{z3'}{$z_{23}$}
\psfrag{z~3}{$ z_{32}$}
\psfrag{z~3'}{$ z_{41}$}
\psfrag{q1}{$p_3$}
\psfrag{q2}{$p_4$}
\psfrag{p1}{$p_1$}
\psfrag{p2}{$p_2$}
{\scalebox{.8}{\includegraphics[height=8cm,width=8cm]{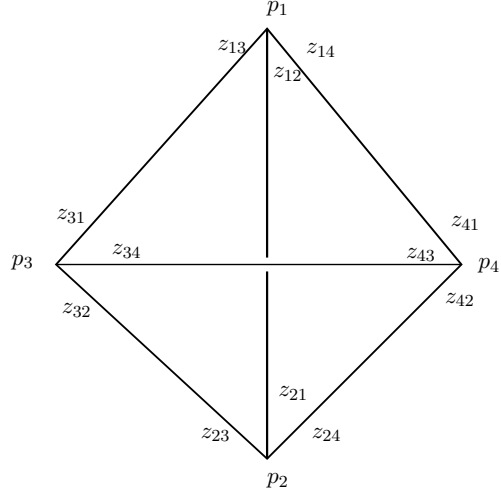}}}
\end{picture}
\end{center}
\caption{\sl Parameters for a cross-ratio structure}\label{Figure:ptetrahedron}
\end{figure}

{\bf Remarks}
\begin{enumerate}
\item To visualize the definition we refer to Figure
 \ref{Figure:ptetrahedron}.
For each $[p_i,p_j,p_k,p_l]\in O\Delta^+$ we define
$$
z_{ij}={\bf X}(p_i,p_j,p_k,p_l).
$$
We interpret $z_{ij}$ as a cross-ratio associated to the edge $[ij]$ at the vertex $i$. 
Cross-ratios of elements of  $O\Delta^-$ are obtained taking inverses by the first 
symmetry.
\item A cross-ratio structure defined on $\Delta$ is a point in the variety in 
$\left ( \C\setminus \{0,1\}\right )^{12}$
with coordinates $z_{ij}$, $1\leq i \neq j\leq 4$ defined by
the usual similarity constraints: if
$(i,j,k,l)$ is an even permutation of $(1,2,3,4)$ then
$$
z_{ik}=\frac{1}{1-z_{ij}}.
$$
\item The similarity relations can be used to reduce the number of variables to four, 
one for each vertex.  One can use, for instance, $(z_{12},z_{21},z_{34},z_{43})\in
\left ( \C\setminus \{0,1\}\right )^{4}$.
\end{enumerate}
In the following we will denote by a sequence 
of numbers $ijkl$ the corresponding invariant ${\bf X}(u_i,u_j,u_k,u_l)$.
Given a simplex $[u_1,u_2,u_3,u_4]$, the twelve coordinates of a cross-ratio structure 
introduced above can be listed as follows
$$\left (\begin{matrix}
z_{12}  \\
z_{13}  \\
z_{14} \\
z_{21}\\
z_{24}\\
z_{23}\\
z_{34}\\
z_{31}\\
z_{32}\\
z_{43}\\
z_{42}\\
z_{41}\\
\end{matrix}\right )=
\left (\begin{matrix}
{1234}  \\
{1342}  \\
{1423} \\
{2143}\\
{2431}\\
{2314}\\
{3412}\\
{3124}\\
{3241}\\
{4321}\\
{4213}\\
{4132}\\
\end{matrix}\right )
$$
By the similarity relations, the cross-ratio structure is defined by 
$z_{12},z_{21},z_{34},z_{43}$ so it is convenient to use the following 
notation.
\begin{dfn}
$$
[[u_1,u_2,u_3,u_4]]=
\left (\begin{matrix}
z_{12}  \\
z_{21}\\
z_{34}\\
z_{43}
\end{matrix}\right )=
\left (\begin{matrix}
{1234}  \\
{2143}\\
{3412}\\
{4321}\\
\end{matrix}\right )
$$
\end{dfn}

\subsection{Hyperbolic configurations}\label{section:hyperbolic}
This section is not used in the next sections.    It identifies 
a subset of cross-ratio structures closely related to real hyperbolic
ideal tetrahedra.  I thank J. Genzmer (see \cite{Ge} for more details)
for correcting an earlier version of it. 

\begin{prop}
The complex algebraic  variety in 
$\left ( \C^*\setminus \{1\}\right )^{12}$
with coordinates $z_{ij}$,  $1\leq i\neq j\leq 4$, defined by, for
$(i,j,k,l)$ an even permutation of $(1,2,3,4)$, 
the usual similarity constraints
$$
z_{ik}=\frac{1}{1-z_{ij}}
$$
and the three complex equations
\begin{eqnarray}\label{eq:hyperbolic}
z_{ij}z_{ji}= {z_{kl}z_{lk}}
\end{eqnarray}
has two irreducible components:
\begin{itemize}
\item One branch is
 parametrised by 
$(z_{12},z_{21})\in \left ( \C^*\setminus \{1\}\right )^{2}$:
$$
z_{34}=- z_{12}\frac{1-z_{21}}{1-z_{12}}
\ \ \ 
z_{43}=- z_{21}\frac{1-z_{12}}{1-z_{21}}.
$$
\item The other branch is parametrised by $\C^*\setminus \{1\}$:
$$
z_{12}=z_{21}=z_{34}=z_{43}.
$$
\end{itemize}
\end{prop}
\Pf By the previous remarks the cross-ratio variety defined by the similarity constraints
 is 
parametrized by $(z_{12},z_{21},z_{34},z_{43})\in
\left ( \C\setminus \{0,1\}\right )^{4}$. It suffices now to solve equations
\ref{eq:hyperbolic} in these coordinates. The equations
are 
$$
z_{12}z_{21}= {z_{34}z_{43}}
$$
$$
(1-z_{12})(1-z_{34})=
(1- z_{21})(1- z_{43})
$$
$$
(1-z_{12})(1-z_{43}){ z_{21}}{z_{32}}=
(1- z_{21})(1- z_{34}){z_{12}}{z_{43}}
$$
and  dividing the last by the second and the first, we obtain
$$
 (1-z_{12})^2{ z_{34}^2}=
(1- z_{21})^2 z_{12}^2
$$
so
$$
z_{34}=\pm z_{12}\frac{1-z_{21}}{1-z_{12}}.
$$
This equation and the first one give then
$$
z_{43}=\pm z_{21}\frac{1-z_{12}}{1-z_{21}}.
$$
Substituting the solution $z_{34}= z_{12}\frac{1-z_{21}}{1-z_{12}}$, 
$z_{43}= z_{21}\frac{1-z_{12}}{1-z_{21}}$
back in the second equation we obtain 
that 
$z_{12}=z_{21}=z_{34}=z_{43}$.
On the other hand the solution  $z_{34}=- z_{12}\frac{1-z_{21}}{1-z_{12}}$, 
$z_{43}=- z_{21}\frac{1-z_{12}}{1-z_{21}}$ satisfies the second equation 
without constraints.
\EPf

Note that the branch 
$$
z_{34}=- z_{12}\frac{1-z_{21}}{1-z_{12}}\ \ \ 
z_{43}=- z_{21}\frac{1-z_{12}}{1-z_{21}}
$$
has the property that
the coeficients
$$
z_1= z_{12}z_{21}, \ z_2= z_{31}z_{13},\ z_3= z_{14}z_{41}
$$
satisfy the similarity conditions
$$
z_{2}=\frac{1}{1-z_1}\ \ {\mbox{and}} \ z_{3}=\frac{1}{1-z_2}.
$$
 There exists, for each hyperbolic configuration
$z_1,z_2,z_3$, a $\C\setminus \{ 0,1\}$-parameter lift in the complex
variety $(\C\setminus \{ 0,1\})^{12}$. 
In particular, the configurations 
$z_1^{1/2}=z_{12}=z_{21}=-z_{34}=-z_{43}$ are in 2-1  correspondence
to ideal hyperbolic configurations.

On the other hand the branch given by $z_{12}=z_{21}=z_{34}=z_{43}$ 
satisfy the condition
$$
z_1z_2z_3=1
$$
which should be opposed to $z_1z_2z_3=-1$ in the hyperbolic case.  Part of these configurations can be interpreted geometrically.  Namely, the real points 
($z_{12}\in \R^*\setminus\{1\}$)
parametrise  configurations of points, up to the action of $PU(2,1)$,  in $S^3\in \C^2$, contained in an $\R$-circle (cf. Lemma 3 in \cite{W1} and section \ref{section:CR}).

\vspace{1cm}

{\bf Remark:} 
Define the map $\left ( \C\setminus \{0,1\}\right )^{12}\rightarrow \left ( \C\setminus \{0\}\right )^{6}$
 by taking $a_{ij}=z_{ij}z_{ji}$. Its image is of complex codimension 2 and
 outside the hyperbolic configurations the map is injective onto 
its image.  On the other hand, the fiber above an ideal hyperbolic configuration is
$\C\setminus \{0,1\}$ as computed in the proposition above. 

\section{Triangulations}
Let $T$ be an ideal triangulation of a 3-manifold. By this we
 mean a simplicial
complex whose underlying topological space is a manifold if 
the vertices are deleted.
Let ${\bf X}(p_i,p_j,p_k,p_l)$ be a function defined on the simplices.

We impose the following compatibility conditions:

\begin{enumerate}
\item Edge compatibility: If $[p_i,p_j,p_{m_0},p_{m_1}], 
[p_i,p_j,p_{m_1},p_{m_2}],\cdots , 
[p_i,p_j,p_{m_n},p_{m_0}]$ are simplices having the  edge $[p_i,p_j]$ 
in common then 
$$
{\bf X}(p_i,p_j,p_{m_0},p_{m_1})\cdots {\bf X}(p_i,p_j,p_{m_n},p_{m_0}) =1
$$
\item Face compatibility:  If $[p_i,p_j,p_k,p_l]$ and $[p_{i'},p_j,p_k,p_l]$
are two simplices with a common face $[p_j,p_k,p_l]$ then
$$
{\bf X}(p_j,p_i,p_k,p_l){\bf X}(p_k,p_i,p_l,p_j){\bf X}(p_l,p_i,p_j,p_k)
=
{\bf X}(p_j,p_{i'},p_k,p_l){\bf X}(p_k,p_{i'},p_l,p_j){\bf X}(p_l,p_{i'},p_j,p_k)
$$
\end{enumerate}

{\bf Remarks}
\begin{enumerate}
\item A more constrained definition of cross-ratio structure 
associated to a triangulation 
 is to give a function  ${\bf X}(p_i,p_j,p_k,p_l)$
defined on all configurations of four distinct vertices in the complex,
 not necessarily in the same simplex. The cross-ratios are then defined 
on the product of four copies of the 0-skeleton.
  The compatibility conditions are more 
difficult to verify in that case but are equivalent to the previous definition
in the case of hyperbolic and CR geometry because the cross-ratios
 define actual
vertices in an appropriate model ( $\C P^1$ and $S^3$ respectively)
 and conversely any four points in the model define a simplex
with cross-ratios defined by these points. One can define an even more constrained condition, which is clearly satisfied for both geometries, by imposing that 
we can add one vertex to the 0-skeleton of the complex keeping the compatibility relations.

In the following section we impose constraints on a triangulation in order
that, once compatibility conditions are verified for a triangulation, they
are valid for other triangulations obtained through Pachner moves.  

\item The definition here should be compared with cross-ratio coordinates
as in \cite{F}.
The relation is that the coordinates, say ${\bf Y}$, defined in  \cite{F} are
given by
$$
{\bf Y}(p_i,p_j,p_k,p_l)={\bf X}(p_i,p_j,p_k,p_l){\bf X}(p_j,p_i,p_l,p_k).
$$ 
The definition of a {\bf T}-structure is essentially the same but the 
computations seem to be more natural. 
 As it will be shown bellow, the
relation can be interpreted as a blow-up of cross-ratio coordinates
along ideal real hyperbolic tetrahedra.

\end{enumerate}

\section{Configurations of five points}

The goal of this section is to obtain all relations between  cross-ratios
obtained from choosing four among five points.  These relations will be used 
to prove a five term relation satisfied by the dilogarithm in the next section.

In order to simplify certain formulae, we 
sometimes denote   ${\bf X}(p_i,p_j,p_k,p_l)$ by $(ijkl)$ or 
simply by $ijkl$.
So, we are looking for all relations between all quadruples $(ijkl)$ with pairwise distinct 
$i,j,k,l$ chosen among ${1,2,3,4,5}$.
We describe the relations by explicitly writing the following formal sum with a minimal
set of independent variables
$$
[[u_1,u_2,u_3,u_4]]-[[u_1,u_2,u_3,u_5]]+[[u_1,u_2,u_4,u_5]]-[[u_1,u_3,u_4,u_5]]+[[u_2,u_3,u_4,u_5]]
$$
$$
=\left (\begin{matrix}
{1234}  \\
{2143}\\
{3412}\\
{4321}\\
\end{matrix}\right )
                          -\left (\begin{matrix}
{1235}  \\
{2153}\\
{3512}\\
{5321}\\
\end{matrix}\right )
                         +\left (\begin{matrix}
{1245}  \\
{2154}\\
{4512}\\
{5421}\\
\end{matrix}\right )
                         -\left (\begin{matrix}
{1345}  \\
{3154}\\
{4513}\\
{5431}\\
\end{matrix}\right )
                           +\left (\begin{matrix}
{2345}  \\
{3254}\\
{4523}\\
{5432}\\
\end{matrix}\right )
$$

Let $\{ i,j,k,l,m\}=\{ 1,2,3,4,5\}$.  

We impose the following relations:
\begin{enumerate}
\item The edge compatibility conditions
$$
(ijkl)=(ijkm)(ijml).
$$
\item The face compatibility conditions 
$$
(ijkl)(ljik)(kjli)=(imkl)(lmik)(kmli).
$$
\end{enumerate}
The edge equations correspond to (oriented) pairs $(ij)$ of points.
That makes 20 equations not all of them independent.
The number of face equations correspond to  (non-oriented) pairs
$\{j,m\}$
 of points
(opposite to each face). That makes 10 equations not all of them independent.

In principle there are 4 independent invariants for each configuration of four
points.  Given five points, there are then $4\times 5=20$ invariants arranged 
in five columns. 

Using the compatibility conditions we may establish relations 
between the columns. In  Proposition \ref{1} we first relate the cross-ratio
invariants using only the edge compatibilities.   
 In fact, the invariants of two sets of configurations
of four points determine those of the other three as is shown in
 Proposition \ref{2} by adding the face compatibilities.  
The proofs are  straightforward computations with the compatibility relations.

\begin{prop} \label{1} If the only relations between the cross-ratios are
the edge compatibilities,  the space of cross-ratios of a configuration of five points
 is of dimension 
10. A possible set of coordinates is 
$$
\{x_1,x_2,x_3,y_1,y_2,y_3,z_3,z_4,w_3,w_4\}:
$$
$$
\left (\begin{array}{c}
x_1  \\
x_2  \\
x_3  \\
\frac{z_3(1-w_3)}{w_3(1-z_3)}  \\
\end{array}\right )
                          -\left (\begin{matrix}
y_1  \\
y_2  \\
y_3  \\
 \frac{1-w_4}{1-z_4} \\
\end{matrix}\right )
                         +\left (\begin{matrix}
\frac{y_1}{x_1}  \\
\frac{y_2}{x_2}  \\
\frac{w_3}{z_3}  \\
\frac{w_4}{z_4}  \\
\end{matrix}\right )
                         -\left (\begin{matrix}
\frac{1-y_1}{1-x_1} \\
\frac{(1-y_3)}{(1-x_3)}\\
w_3\\
w_4\\
\end{matrix}\right )
                           +\left (\begin{matrix}
\frac{x_2(1-y_2)}{y_2(1-x_2)} \\
\frac{x_3(1-y_3)}{y_3(1-x_3)}\\
z_3\\
z_4\\
\end{matrix}\right )
$$
\end{prop}
\Pf  The proof follows writing all edge compatibility relations.  For instance,
the 4th element in the first column is $(4321)$.  We have
$$
(4321)= (4325)(4351)=\frac{1}{1-\frac{1}{z_3}}(1-\frac{1}{w_3})=\frac{z_3(1-w_3)}{w_3(1-z_3)}.
$$
Also, the 4th element in the second column is $(5321)$.  We have
$$
(5321)= (5324)(5341)=\frac{1}{1-z_4}(1-w_4).
$$
The other terms are obtained similarly.
\EPf

\begin{prop}\label{2} If the only relations between the cross-ratios  are
 edge and face compatibilities,  the space of cross-ratios of a configuration of five points is of dimension 7. 
A possible set of coordinates is 
$$
\{x_1,x_2,x_3,x_4,y_1,y_2,y_4\}:
$$
$$
\left (\begin{array}{c}
x_1  \\
x_2  \\
x_3  \\
x_4  \\
\end{array}\right )
                          -\left (\begin{matrix}
y_1  \\
y_2  \\
y_3  \\
y_4  \\
\end{matrix}\right )
                         +\left (\begin{matrix}
\frac{y_1}{x_1}  \\
\frac{y_2}{x_2}  \\
\frac{w_3}{z_3}  \\
\frac{w_4}{z_4}  \\
\end{matrix}\right )
                         -\left (\begin{matrix}
\frac{1-y_1}{1-x_1} \\
\frac{(1-y_3)}{(1-x_3)}\\
w_3\\
w_4\\
\end{matrix}\right )
                           +\left (\begin{matrix}
\frac{x_2(1-y_2)}{y_2(1-x_2)} \\
\frac{x_3(1-y_3)}{y_3(1-x_3)}\\
z_3\\
z_4\\
\end{matrix}\right )
$$
with
$$
y_3=\frac{y_1x_3(y_2-1)(x_1-1)}{x_1(y_1-1)(x_2-1)}.
$$
$$
z_3=\frac{x_4(-y_1x_1y_2-x_1x_2+y_1x_1x_2+x_1+y_1y_2-y_1)}
{y_1(x_4-1)(x_2-y_2)(x_1-1)}
$$
$$
w_3=
\frac{-y_1x_1y_2-x_1x_2+y_1x_1x_2+x_1+y_1y_2-y_1}
{(x_4-1)(x_1-y_1)(x_2-1)}
$$
$$
z_4=\frac{x_1(y_4-1)(y_1-1)(x_2-y_2)}
{y_4(-y_1x_1y_2-x_1x_2+y_1x_1x_2+x_1+y_1y_2-y_1)}.
$$
$$
w_4=\frac{(y_4-1)(y_2-1)(x_1-y_1)}{-y_1x_1y_2-x_1x_2+y_1x_1x_2+x_1+y_1y_2-y_1}.
$$
\end{prop}\label{prop:5term}
\Pf
It remains to use the face compatibilities:
\begin{enumerate}
\item
From 
$$
(1524)(4512)(2541)=(1324)(4312)(2341)
$$ we get
$$
\frac{(x_1-y_1)x_2}{y_1(z_3x_4-x_4-z_3)(-y_2+x_2)}=
-\frac{x_2(-1+x_1)}{(x_2-1)x_4}
$$
Therefore
$$
z_3=\frac{x_4(-y_1x_1y_2-x_1x_2+y_1x_1x_2+x_1+y_1y_2-y_1)}
{y_1(x_4-1)(x_2-y_2)(x_1-1)}.
$$
Substituting the above value for $z_3$ in  $x_4=\frac{z_3(1-w_3)}{w_3(1-z_3)}$
we obtain 
$$
w_3=
\frac{-y_1x_1y_2-x_1x_2+y_1x_1x_2+x_1+y_1y_2-y_1}
{(x_4-1)(x_1-y_1)(x_2-1)}
$$
\item Analogously, from
$$
(1452)(2415)(5421)=(1352)(2315)(5321)
$$ we get
$$
z_4=\frac{x_1(y_4-1)(y_1-1)(x_2-y_2)}
{y_4(-y_1x_1y_2-x_1x_2+y_1x_1x_2+x_1+y_1y_2-y_1)}.
$$
Again, substituting the above value for $z_4$ in $y_4=\frac{1-w_4}{1-z_4}$ we
obtain
$$
w_4=\frac{(y_4-1)(y_2-1)(x_1-y_1)}{-y_1x_1y_2-x_1x_2+y_1x_1x_2+x_1+y_1y_2-y_1}.
$$

\item The equation
$$
(1523)(3512)(2531)=(1423)(3412)(2431)
$$
gives
$$
y_3=\frac{y_1(y_2-1)(x_1-1)x_3}{(y_1-1)(x_2-1)x_1}.
$$
\item A computation shows that 
the other equations don't give any new relations.
\end{enumerate}
\EPf

\section{Dilogarithm and volume}
In this section we define a volume of a cross-ratio structure on a simplex.
For preliminaries on the dilogarithm we refer the reader to \cite{Z}. Consider the function (Bloch-Wigner) (see \cite{B}, section 3 of \cite{Z} or formula 19 in \cite{O})
$$
D(z)= log|z| arg(1-z) -\Im \int_0^z\frac{log(1-t)}{t} dt
$$
which is well defined and analytic on $\C\setminus\{0,1\}$ and
extends to a continuous function on $\C P^1$ by defining 
$D(0)=D(1)=D(\infty)=0$.
It satisfies the 5-term relation (see formula 34 in \cite{O})
$$
D(x)-D(y)+D(\frac{y}{x})-D(\frac{1-y}{1-x})+D(\frac{1-y^{-1}}{1-x^{-1}})=0.
$$
There are many equivalent forms of the five term relations. Each one is 
obtained from the other
by a change of coordinates. For instance the five term relation in \cite{Z}, formula 4 is
$$
D(u)+D(v)+D(\frac{1-u}{1-uv})+D({1-uv})+D(\frac{1-v}{1-uv})=0
$$
which can be obtained from the previous by writing $v=1/y$.

Recall that the Bloch-Wigner function can be interpreted as a volume function on the
space of ideal hyperbolic tetrahedra (see section 4 in \cite{Z}).  Indeed, 
$\C\setminus\{0,1\}$ parametrises configurations of four distinct points in $\C P^1$
which is identified to the boundary of real hyperbolic space, $H^3_\R$.  The convex hull 
(inside $H^3_\R$) of four points in $\C P^1$  with cross-ratio $z$  defines an ideal simplex
, up to translations by $PSL(2,\C)$, whose volume is $D(z)$.

We will define next a function defined on cross-ratio structures.

Associated to a cross-ratio structure are the invariants $z_{ij}$.
Recall that four invariants, one  at each vertex, determine the whole set of invariants, so 
the we might chose  
$z_1=z_{12}$, $z_2=z_{21}$, $z_3=z_{34}$, $z_4=z_{43}$.  It is reasonable to expect
that the   following
definition will be an analog of the volume of an ideal hyperbolic simplex.  But the true
reason behind it will be the fact that it satisfies a 5 term relation.

\begin{dfn}\label{dfn:volume}
The volume of a cross-ratio structure ${\bf z}=(z_1,z_2,z_3,z_4)\in \left ( \C\setminus\{0,1\}\right ) ^4 $ on a simplex $[p_1,p_2,p_3,p_4]$ is
$$
{\cal D}({\bf z})= D(z_1)+D(z_2)+D(z_3)+D(z_4).
$$
\end{dfn}

Using  formula (see \cite{Z}, formula 2)
$$
D(z)=\frac{1}{2}\left (D\left(\frac{z}{\bar z}\right)+D\left(\frac{1-z^{-1}}{1-\bar z^{-1}}\right)
+D\left(\frac{1-\bar z}{1-z}\right)\right )
$$
we obtain (cf. Lemma 2 for  hyperbolic geometry in \cite{M}) 
$$
{\cal D}({\bf z})=\frac{1}{2}\sum_{ij}D(e^{2i \theta_{ij}})=\sum_{ij}\Lambda(\theta_{ij}),
$$
where $\theta_{ij}=\arg\, z_{ij}$ and $\Lambda(\theta)$ is Lobachevsky function as defined
by \cite{Co,M}.
\vspace{1cm}

Let $T$ be an ideal triangulation with a cross-ratio structure as above satisfying
edge and face compatibilities.
We define a function  on simplices of the triangulation
 by using the volume function of a generic tetrahedron.
Observe that to each simplex $[p_1,p_2,p_3,p_4]$ we associate four complex
 coordinates ${\bf z}=(x_1,x_2,x_3,x_4)\in \left ( \C\setminus\{0,1\}\right ) ^4$.
If the cross-ratio structure is fixed we will also write 
$$
{\cal D}([p_1,p_2,p_3,p_4])= {\cal D}({\bf z})=D(x_1)+D(x_2)+D(x_3)+D(x_4).
$$
\begin{thm}\label{theorem:D}
The function ${\cal D}$ satisfies the 5-term relation:
$$
{\cal D}([p_1,p_2,p_3,p_4])
-{\cal D}([p_1,p_2,p_3,p_5])
+{\cal D}([p_1,p_2,p_3,p_4])
-{\cal D}([p_1,p_3,p_4,p_5])
+{\cal D}([p_2,p_3,p_4,p_5])=0.
$$
\end{thm}
\Pf
Recall that
$$
D(x)-D(y)+D(\frac{y}{x})-D(\frac{1-y}{1-x})+D(\frac{1-y^{-1}}{1-x^{-1}})=0.
$$
It gives rise 
 to the following relations:
$$
D(x_1)-D(y_1)+D(\frac{y_1}{x_1})-D(\frac{1-y_1}{1-x_1})+
D(\frac{1-y_1^{-1}}{1-x_1^{-1}})=0,
$$
$$
D(x_2)-D(y_2)+D(\frac{y_2}{x_2})-D(\frac{1-y_2}{1-x_2})+
D(\frac{1-y_2^{-1}}{1-x_2^{-1}})=0,
$$
$$
D(x_3)-D(y_3)+D(\frac{y_3}{x_3})-D(\frac{1-y_3}{1-x_4})+
D(\frac{1-y_3^{-1}}{1-x_3^{-1}})=0,
$$
$$
D(x_4)-D(y_4)+D(\frac{y_4}{x_4})-D(\frac{1-y_4}{1-x_4})+
D(\frac{1-y_4^{-1}}{1-x_4^{-1}})=0.
$$
We also have
$$
D(\frac{1-w_3^{-1}}{1-z_3^{-1}})-D(\frac{1-w_3}{1-z_3})+
D(\frac{w_3}{z_3})- D(w_3)+D(z_3)=0
$$
and
$$
D(\frac{1-w_4^{-1}}{1-z_4^{-1}})-D(\frac{1-w_4}{1-z_4})+
D(\frac{w_4}{z_4})- D(w_4)+D(z_4)=0.
$$
Comparing these relations to the sum of  the lines
in the formal sum of  Proposition \ref{prop:5term} we obtain
$$
{\cal D}([p_1,p_2,p_3,p_4])
-{\cal D}([p_1,p_2,p_3,p_5])
+{\cal D}([p_1,p_2,p_3,p_4])
-{\cal D}([p_1,p_3,p_4,p_5])
+{\cal D}([p_2,p_3,p_4,p_5])
=
$$
$$
-D(\frac{1-y_1^{-1}}{1-x_1^{-1}})
+D(\frac{1-y_2}{1-x_2})
-D(\frac{y_3}{x_3})
+D(\frac{1-w_3}{1-z_3})
-D(\frac{1-w_4^{-1}}{1-z_4^{-1}})
$$
Now, substituting the values of $y_3,z_3, w_3,z_4, w_4$ in terms 
of $x_1,x_2,x_3,x_4,y_1,y_2,y_4$ we obtain
$$
-D(\frac{1-y_1^{-1}}{1-x_1^{-1}})
+D(\frac{1-y_2}{1-x_2}) 
-D(\frac{(-1+x_1)(-y_2+x_2)y_1}{(x_1-y_1)(x_2-1)})
+D(\frac{x_1(-1+y_1)(-y_2+x_2)}{(x_1-y_1)(y_2-1)}).
$$
Call $a=\frac{1-y_1^{-1}}{1-x_1^{-1}}$ and $b=\frac{1-y_2}{1-x_2}$. Then,
a simple calculation gives that the above expression is
$$
-(D(a)-D(b)+D(\frac{b}{a})-D(\frac{1-b}{1-a})+D(\frac{1-b^{-1}}{1-a^{-1}}))=0.
$$
\EPf

Consider a triangulation of a three manifold with a cross-ratio structure.
The volume is obtained by adding over
all simplices:
$$
Vol=\sum_i \epsilon_i {\cal D}(T_i)
$$ 
where the sum is taken over all 3-simplices with a factor $\epsilon=\pm 1$
which is $+1$ if the orientation of the simplex is the same as 
of the space and  $-1$ if the orientation of the simplex is the  opposite. 
 By the five term relation, it does not depend on the 
triangulation.
The volume defines an element in $H^3(M,\R)$.  See \cite{T,NZ} for the case of real hyperbolic geometry.
\vspace{.5cm}

{\bf Remark}:
For the hyperbolic configurations, that is, when 
$$
z_{34}=- z_{12}\frac{1-z_{21}}{1-z_{12}}\ \ \ 
z_{43}=- z_{21}\frac{1-z_{12}}{1-z_{21}},
$$
a simple use of the five term relation shows that
$$
D(z_{12}z_{21})=D(z_{12})+D(z_{21})+D(z_{34})+D(z_{43}).
$$
This shows that the function defined above coincides with the usual
volume function for ideal hyperbolic tetrahedra.

\section{CR geometry (see \cite{BS,G,J})}\label{section:CR}

\label{CRgeometry}

CR geometry is modeled on the {\sl Heisenberg group} $\Heis$,
 the set of pairs $(z, t)\in
{\C}\times{\R}$ with the product 
$$
(z,t)\cdot (z',t') = (z+z', t + t' + 2 \Im z \overline{z}').
$$ 
The one point compactification of the Heisenberg group, 
$\overline{\Heis}$, of $\Heis$ can be interpreted as $S^3$ which
, in turn, can be identified to the boundary of Complex Hyperbolic
space.   

We consider the group $U(2,1)$ preserving the Hermitian form 
$\langle z,w \rangle = w^*Jz$ defined by 
the matrix
$$
J=\left ( \begin{array}{ccc}

                        0      &  0    &       1 \\

                        0       &  1    &       0\\

                        1       &  0    &       0

                \end{array} \right )
$$ 
and the following subspaces in ${\mathbb C}^{3}$:
$$
        V_0 = \left\{ z\in {\mathbb C}^{3}-\{0\}\ \ :\ \
 \langle z,z\rangle = 0 \ \right\},
$$
$$
        V_-   = \left\{ z\in {\mathbb C}^{3}\ \ :\ \ \langle z,z\rangle < 0 
\ \right\}.
$$
Let $\PP:{\C}^{3}\setminus\{ 0\} \rightarrow {\C}P^{2}$ be the
canonical projection.  Then
${\bf H}_{\C}^{2} = \PP(V_-)$ is the complex hyperbolic space and 
$S^3={\bf H}_{\C}^{2} = \PP(V_0)$ can be identified to $\overline\Heis$.

The group of biholomorphic transformations of ${\bf H}_{\C}^{2}$ is then 
$PU(2,1)$, the projectivization of $U(2,1)$.  It acts on $S^3$ by 
CR transformations. 
  We 
define $\C$-circles as boundaries of complex lines in ${\bf H}_{\C}^{2}$. Analogously,
$\R$-circles are boundaries of totally real totally geodesic two dimensional submanifolds in ${\bf H}_{\C}^{2}$.
Using the identification $S^3= \Heis \cup \{ \infty\}$
one can define alternatively a $\C$-circle
 as any circle in $S^3$ which is obtained
from the vertical line $\{(0,t)\}\cup \{ \infty\}$ in the compactified Heisenberg
space by translation by an element of $PU(2,1)$.  Analogously, $\R$-circles are all obtained
by translations of the horizontal line $\{(x,0)\}\cup \{ \infty\}$, $x\in \R$.

A point $p=(z,t)$ in the Heisenberg group and the point $\infty$ are lifted 
to the following points in $\C^{2,1}$:
$$
\hat{p}=\left[\begin{matrix}\frac{-|z|^2+it}{2} \\ z \\ 1 \end{matrix}\right]
\quad\hbox{ and }\quad
\hat{\infty}=\left[\begin{matrix} 1 \\ 0 \\ 0 \end{matrix}\right].
$$
\begin{dfn}
Given any three ordered points $p_1$, $p_2$, $p_3$
 in $\partial\ch{2}$ we define 
{\sl Cartan's angular invariant} $\A$ as 
$$
\A(p_1,p_2,p_3)=
\arg(-\langle\hat{p}_1,\hat{p}_2\rangle\langle\hat{p}_2,\hat{p}_3\rangle
\langle\hat{p}_3,\hat{p}_1\rangle).
$$
\end{dfn}
The Cartan's angular invariants classifies ordered triples of points
in $S^3$:
\begin{prop}[\cite{C}, see also\cite{G}]
There exists an element of $PU(2,1)$ which translates an ordered triple
of points in $S^3$ to another if and only if their corresponding Cartan's
invariants are equal.
\end{prop}

The CR cross ratio is given by
the Koranyi-Reimann invariant introduced in \cite{KR} (see \cite{KR} and
\cite{G} for its properties):
\begin{dfn} The CR cross-ratio associated to four distinct points in $S^3$
is
$$
KR(p_1,p_2,p_3,p_4)=\frac{\langle p_4,p_2\rangle \langle p_3,p_1\rangle}
{\langle p_3,p_2\rangle \langle p_4,p_1\rangle}.
$$
\end{dfn}
Here, we choose lifts for the points $p_i$ which we denote 
by the same letter.  The invariant does not
depend on the choice of lifts. 
The product of the three cross ratios gives the Cartan invariant
(see \cite{KR,G})
$$
KR(p_1,p_2,p_3,p_4)KR(p_1,p_4,p_2,p_3)KR(p_1,p_3,p_4,p_2)= 
e^{2i\A (p_2,p_3,p_4)}
$$

\subsection{Configurations of four points}
We refer to Figure \ref{Figure:ptetrahedron} to describe the
parameters of a tetrahedron (see also \cite{F}). 
Consider a generic configuration
of four (ordered) points in $S^3$ (any three of them
not contained in a $\C$-circle).
Fix one of them 
say $p_1$ and consider the 
projective space of complex lines passing through it.  Then $p_2,p_3,p_4$
 determine three points $t_2, t_3,t_4$ on $\C P^1$.  
The fourth point corresponds to the
complex line passing through $p_1$ and tangent to $S^3$, call it $t_1$.  The 
cross-ratio
of those four points in $\C P^1$ is $z_{12}={\bf X}(t_1,t_2,t_3,t_4)$ (here, ${\bf X}$
is the usual cross-ratio of four points in $\C P^1$).  We define analogously 
the other invariants.  If we take $p_1=\infty$, the complex lines passing through $p_1$
intersect $\Heis$ in vertical lines which are then determined by a coordinate in $\C$.
Up to Heisenberg translations, we can assume that $p_2=(0,0)$ and $p_3=(1,s_3)$
and $p_2=(z_{12},s_4)$, $s_3,s_4\in \R$. The corresponding points in $\C P^1$ will be
$\infty,0,1,z_{12}$.  Therefore one ``sees''
at the vertex $p_1$ the  Euclidean triangle determined by
$0,1,z_{12} \in \C$.

We associate to each vertex $i\in [ij]$ inside an edge the invariant
$(ijkl)$ where the order $k$ and $l$ is fixed by the right hand with the thumb
pointed from $j$ to $i$. A shortcut notation for the invariants is
therefore
$$
z_{ij}=(ijkl),
$$
the indices $kl$ being determined by the choice $ij$.
It satisfies the relation $(ijlk)=(ijkl)^{-1}$.

They satisfy the following relations:
$$
z_{ij}z_{ji}=\overline {z_{kl}z_{lk}} 
$$

{\bf Remarks:} 
\begin{enumerate}
\item An explicit formula for the invariants is given in \cite{W2};
$$
(p_1,p_2,p_3,p_4)=\frac{\langle \hat p_4, c_{12}\rangle \langle \hat p_3, \hat p_1 \rangle}
{\langle \hat p_3, c_{12}\rangle \langle \hat p_4, \hat p_1 \rangle},
$$
where $\hat p_i$ are lifts of   $p_i$ and $c_{12}\in \C^{2,1}$
is a vector orthogonal to the complex plane defined by   $\hat p_1$ and  $\hat p_2$.
\item  The formulae
 relating the cross-ratio invariants
and those defined above are
$$
(p_1,p_2,p_3,p_4)=\frac{ KR(p_1,p_2,p_3,p_4)KR(p_1,p_3,p_4,p_2)KR(p_2,p_3,p_1,p_4)+1}{1+KR(p_1,p_4,p_2,p_3)(KR(p_4,p_2,p_1,p_3)-1)}
$$
and conversely
$$
KR(p_1,p_2,p_3,p_4)=(p_1,p_2,p_3,p_4)(p_2,p_1,p_4,p_3).
$$
\item For other descriptions of  configurations of four points
in $S^3$  and their applications we refer to \cite{W1, W2, PP, PP1,FP}.
\end{enumerate}

\begin{prop} (cf. \cite{F})
Configurations of four distinct points in $S^3$ such that any three points are not contained 
in a $\C$-circle or an $\R$-circle
are parametrised by 
the real algebraic  variety in 
$\left ( \C ^*\setminus \{ 1\}\right )^{12}\setminus \left( \R^*\setminus \{ 1\}\right )^{12}$
with coordinates $z_{ij}$,  $1\leq i\neq j\leq 4$, defined by, for
$(i,j,k,l)$ an even permutation of $(1,2,3,4)$, 
the usual similarity constraints
$$
z_{ik}=\frac{1}{1-z_{ij}}
$$
and the three complex equations
\begin{eqnarray}\label{eq:cr}
z_{ij}z_{ji}=\overline {z_{kl}z_{lk}}
\end{eqnarray}
\end{prop}
{\bf Remarks}:
\begin{enumerate}
\item The  real solutions are contained in two different 
branches (cf. Proposition \ref{section:hyperbolic}) .  One parametrises
 configurations with 
 four points 
 contained in an $\R$-circle.  The other branch corresponds to 
 degenerate hyperbolic ideal tetrahedra with  four points contained 
in the boundary of a totally geodesic plane in real hyperbolic space.
 I thank J. Genzmer for correcting an earlier version 
which appeared in \cite{F}.  For more details see \cite{Ge}.
\item 
 In fact, eliminating two variables at each vertex,
 one can write the six real equations directly in 
$\left ( \C\setminus \{0,1\}\right )^{4}$ with variables
$z_{12},z_{21},z_{34},z_{43}$.  The equations in these variables
are:
$$
z_{12}z_{21}= \overline{z_{34}z_{43}}
$$
$$
\frac{1}{1-z_{12}}\frac{1}{1-z_{34}}=
\frac{1}{1-\bar z_{21}}\frac{1}{1-\bar z_{43}}
$$
$$
(1-\frac{1}{z_{12}})(1-\frac{1}{z_{43}})=
(1-\frac{1}{\bar z_{34}})(1-\frac{1}{\bar z_{21}})
$$

\item
Recall that (see \cite{KR,G})
$$
KR(p_1,p_2,p_3,p_4)KR(p_1,p_4,p_2,p_3)KR(p_1,p_3,p_4,p_2)= 
e^{2i\A (p_2,p_3,p_4)}
$$
gives Cartan's invariant in terms of cross ratios.
One can write then
$$
e^{2i\A (p_2,p_3,p_4)}=z_{12}z_{21}z_{14}z_{41}z_{13}z_{31}=-z_{21}z_{41}z_{31}.
$$

A common face of two tetrahedra has opposite orientations and as $\A (p_3,p_2,p_4)=-\A (p_2,p_3,p_4)$,
 the face gluing conditions between tetrahedra with invariants
 $z_{ij}$ and $w_{i'j'}$ are
given by expressions of the form
$$
z_{il}z_{jl}z_{kl}w_{i'l'}w_{j'l'}w_{k'l'}=1.
$$
where $l$ and $l'$ correspond to points oposed to the common face.
This explains the face compatibility conditions in the CR case.
\item
Writing
$$
z_{ij}=r_{ij}e^{i\theta_{ij}}
$$
we observe that the angles $\theta_{ij}$ determine the parameters $z_{ij}$.
The equations defining the possible values of $\theta_{ij}$ are:
\begin{enumerate}
\item For each vertex $i$:
$$
\sum_j \theta_{ij}=\pm \pi.
$$
\item Two sets of three CR conditions (there are only four independent equations, two from each set):
$$
 \theta_{ij}+ \theta_{ji}+ \theta_{kl}+ \theta_{lk} = 0 \ (2\pi)
$$
$$
r_{ij}r_{ji}=r_{kl}r_{lk}.
$$
\end{enumerate}
Using the relations at each vertex of the form $r_{12}=\frac{\sin \theta_{13}}{\sin \theta_{14}}, r_{21}=\frac{\sin \theta_{24}}{\sin \theta_{23}}, \cdots$, 
we may write the last conditions in terms of angles as
$$
\frac{\sin \theta_{13}}{\sin \theta_{14}}\frac{\sin \theta_{24}}{\sin \theta_{23}}=
\frac{\sin \theta_{31}}{\sin \theta_{32}}\frac{\sin \theta_{42}}{\sin \theta_{41}}
$$
$$
\frac{\sin \theta_{14}}{\sin \theta_{12}}\frac{\sin \theta_{32}}{\sin \theta_{34}}=
\frac{\sin \theta_{23}}{\sin \theta_{21}}\frac{\sin \theta_{41}}{\sin \theta_{43}}
$$
$$
\frac{\sin \theta_{12}}{\sin \theta_{13}}\frac{\sin \theta_{43}}{\sin \theta_{42}}=
\frac{\sin \theta_{21}}{\sin \theta_{24}}\frac{\sin \theta_{34}}{\sin \theta_{31}}.
$$
The last equation is clearly obtained from the first two.  
There are 12 variables $\theta_{ij}$, 4 equations at each vertex and 4 equations
corresponding to the CR conditions.  That makes a total of 4 independent parameters.
 \end{enumerate}

\vspace{1cm}

\subsection{The CR volume as a coboundary}

As it was pointed out to me by Qingxue Wang the fact that $H_{cont}^3(PU(2,1),\R)=0$
 implies that ${\cal D}$ (which can be seen as a measurable 3-cocycle in $PU(2,1)$)
 is a coboundary and, therefore,
 the volume function of a CR structure of any closed three manifold is null.

For background on continous cohomology we refer to Lecture 3 in \cite{B1} (a more comprehensive introduction is \cite{Gu}).
The vanishing of the continuous cohomology group follows from Van Est theorem 
that $H_{cont}^3(PU(2,1),\R)= H^3(\mathfrak{ g}, \mathfrak{ u}, \R)$, where $\mathfrak{ g}$
is the Lie algebra of $PU(2,1)$ and $\mathfrak{u}$ the Lie algebra of the maximal compact 
subgroup $U(2)$. Indeed, 
 write a Cartan decomposition $\mathfrak{g}=\mathfrak{ u}+\mathfrak{ p}$
and let $\mathfrak{ g_u}=\mathfrak{ u}+i\mathfrak{ p}$ be another compact form in the complexified Lie algebra. We obtain  (cf. Lecture 3 in \cite{B1} or chapter III, 7 in \cite{Gu}) that
 $H^3(\mathfrak{ g}, \mathfrak{ u}, \R)= H^3(\mathfrak{ g_u}, \mathfrak{ u}, \R)$ and  $H^3(\mathfrak{ g_u}, \mathfrak{ u}, \R)=H^3(\C P^2, \R)=0$ ($\C P^2$ being the compact symmetric
space associated to non-compact symmetric space $H^2_\C$). 

Continuous cohomology can also be computed using measurable cochains (see \cite{B1}).
Fix a point $\infty\in S^3$ and consider the measurable cochain defined
 outside a set of zero measure in $PU(2,1)^4$, equipped with a Haar measure,
by
$$
{\cal D}(g_1,g_2,g_3,g_4)={\cal D}([g_1\infty,g_2\infty,g_3\infty,g_4\infty]).
$$
The set of measure zero, where the 2-cochain is not defined, is the set 
of quadruples such that the points $g_1\infty,g_2\infty,g_3\infty,g_4\infty$
are either not pairwise distinct or degenerate (three of them belong to
a $\C$-circle).  Theorem \ref{theorem:D} is the statement that the measurable cochain
is a measurable cocycle.

In the following we will determine
${\cal D}$ as a coboundary. 
Define the measurable 2-cochain in $PU(2,1)$ by
$$
c_2(g_1,g_2,g_3)=\frac{1}{2}D(-e^{2i\A (g_1\infty,g_2\infty,g_3\infty)}).
$$

{\bf{Remark}}:
We can also define ${\cal D}$ and $c_2$ as cochains for the 
simplicial complex defined by a triangulation of a three manifold
with a cross-ratio structure.  
In that case, if $[p_1,p_2,p_3,p_4]$ is a simplex, ${\cal D}(p_1,p_2,p_3,p_4)=
{\cal D}({\bf z})$ as in Definition \ref{dfn:volume} and
$c_2(p_2,p_3,p_4)=D(z_{21}z_{31}z_{41})/2$.  But in the general case, the volume is not 
a coboundary.  The following proposition is, therefore, special to the CR case.

\begin{prop}
$$
\partial c_2 ={\cal D}
$$
\end{prop}
\Pf Using the definition of a coboundary, we have to prove that
$$
D(z_{21}z_{41}z_{31})+D(z_{12}z_{32}z_{42})+D(z_{13}z_{23}z_{43})+D(z_{14}z_{24}z_{34})=2({ D}(z_{12})+{ D}(z_{21})+{ D}(z_{34})+{ D}(z_{43})).
$$

\begin{figure}
\setlength{\unitlength}{1cm}
\begin{center}
\begin{picture}(7,10)
\psfrag{z1}{$\infty$}
\psfrag{z1'}{$0$}
\psfrag{z~1}{$z$}
\psfrag{z~1'}{$1$}
\psfrag{z2}{$1$}
\psfrag{z2'}{$\infty$}
\psfrag{z~2}{$ \frac{1}{1-z}$}
\psfrag{z~2'}{$ 0$}
\psfrag{z3}{$0$}
\psfrag{z3'}{$1$}
\psfrag{z~3}{$1-\frac{1}{z}$}
\psfrag{z~3'}{$\infty$}
\psfrag{q1}{$p_3$}
\psfrag{q2}{$p_4$}
\psfrag{p1}{$p_1$}
\psfrag{p2}{$p_2$}
{\scalebox{.8}{\includegraphics[height=8cm,width=8cm]{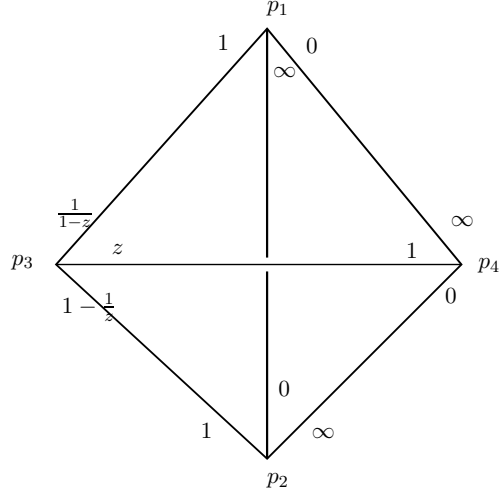}}}
\end{picture}
\end{center}
\caption{\sl Parameters for a degenerate tetrahedron}\label{Figure:degenerate}
\end{figure}

${\cal D}$ is a continuous  3-chain defined on generic tetrahedra which extends
to degenerate tetrahedra of the form of Figure \ref{Figure:degenerate}.
In that case ${\cal D}=D(z)$. Suppose ${\cal D}$ is a coboundary. 
We may suppose that there exists a function $F: U(1)\rightarrow \R$ such that  
$$
F(z_{21}z_{41}z_{31})+F(z_{12}z_{32}z_{42})+F(z_{13}z_{23}z_{43})+F(z_{14}z_{24}z_{34})=2({ D}(z_{12})+{ D}(z_{21})+{ D}(z_{34})+{ D}(z_{43})).
$$
That function might be extended to degenerate tetrahedra and in that case
we compute, taking limits using the CR constraints,
 $$
z_{21}z_{41}z_{31}=\frac{1-\bar z}{1-z}
$$
$$
z_{12}z_{32}z_{42}= \frac{1-\frac{1}{z}}{1-\frac{1}{\bar z}}
$$
$$
z_{14}z_{24}z_{34}=\frac{z}{\bar z}
$$
$$
z_{13}z_{23}z_{43}=1.
$$
Therefore, fixing $F(1)=0$, we have 
$$
F\left(\frac{z}{\bar z}\right)+F\left(\frac{1-\bar z}{1-z}\right)+F\left( \frac{1-\frac{1}{z}}{1-\frac{1}{\bar z}}\right)=2 D(z).
$$
Using angle variables we obtain
$$
F(e^{2i\theta_1})+F(e^{2i\theta_2})+F(e^{2i\theta_3})=
 D(e^{2i\theta_1})+D(e^{2i\theta_2})+D(e^{2i\theta_3})
$$
where $\theta_1+\theta_2+\theta_3 = 0\ (\pi)$.  Taking derivatives with respect
to $\theta_1$ and supposing $\theta_2$ independent we get
$$
(F-D)'(e^{2i\theta_1})-(F-D)'(e^{2i\theta_3})=0.
$$
As $\theta_2$ is independent we conclude that $(F-D)'$ is constant.  In fact,
it is null as the functions are periodic.
This implies that $F= D$ up to an additive constant which must be zero as $F(1)=D(1)=0$.

\EPf

\end{document}